\documentclass[a4paper,10pt,reqno]{amsart}

\usepackage{amssymb,amsmath,amsfonts} 
\usepackage{a4wide}
\usepackage{graphicx}
\usepackage{hyperref}
\usepackage{color}
\usepackage{amsthm}
\usepackage{enumerate}
\numberwithin{equation}{section}    
\theoremstyle{plain}
\newtheorem{Theorem}{Theorem}[section]
\newtheorem{Proposition}[Theorem]{Proposition}
\newtheorem{Corollary}[Theorem]{Corollary}
\newtheorem{Lemma}[Theorem]{Lemma}
\theoremstyle{definition}

\newtheorem{Definition}[Theorem]{Definition}
\newtheorem{Example}[Theorem]{Example}
\newtheorem{Examples}[Theorem]{Examples}

\theoremstyle{remark}
\newtheorem{Remark}[Theorem]{Remark}
\renewcommand{\phi}{\varphi}
\newcommand{\one}{\mathbf{1}} 

\newcommand{\RR}{\mathbb{R}}

\newcommand{\NN}{\mathbb{N}}

\newcommand{\ZZ}{\mathbb{Z}}

\newcommand{\cC}{\mathcal{C}}

\newcommand{\cP}{\mathcal{P}}

\newcommand{\cH}{\mathcal{H}}

\newcommand{\cA}{\mathcal{A}}

\newcommand{\cG}{\mathcal{G}}
\newcommand{\cL}{\mathcal{L}}

\newcommand{\cN}{\mathcal{N}}
\newcommand{\cO}{\mathcal{O}}
\newcommand{\cQ}{\mathcal{Q}}

\newcommand{\hm}[1]{\textbf{*}\leavevmode{\marginpar{\tiny%
$\hbox to 0mm{\hspace*{-0.5mm}$\leftarrow$\hss}%
\vcenter{\vrule depth 0.1mm height 0.1mm width \the\marginparwidth}%
\hbox to 0mm{\hss$\rightarrow$\hspace*{-0.5mm}}$\\\relax\raggedright #1}}}
\title[Convergence theorems for graph sequences]{Convergence theorems for graph sequences}
\author[Felix Pogorzelski]{Felix Pogorzelski}

\begin{document}

\begin{abstract}
In this paper, we deal with a notion of Banach space-valued mappings defined on a set consisting of finite graphs with uniformly bounded vertex degree. These functions will be endowed with certain boundedness and additivity criteria. We examine their normalized long-term behaviour along a particular class of graph sequences. Using techniques developed by Elek, we show convergence in the topology of the Banach space if the corresponding graph sequence possesses a hyperfinite structure. These considerations extend and complement the corresponding results for amenable groups. 
As an application, we verify the uniform approximation of the integrated density of states for bounded, finite range operators on discrete structures. Further, we extend results concerning an abstract version of Fekete's Lemma for amenable groups and cancellative semigroups to the geometric situation of convergent graph sequences. 
\end{abstract}

\maketitle

\begin{center}
{\bf Address} \\
Technion, Israel Institute of Technology \\
Technion Campus, 32000 Haifa, Israel\\
{\tt felixp@technion.ac.il}\\
\end{center}

{\bf Keywords:} Graph convergence; hyperfiniteness; integrated density of states. \\

{Mathematics subject classification:} \quad {05C25, 05C80, 46L60, 47A35}.

\section{Introduction}

The following work is devoted to convergence results for Banach space valued functions with particular (sub-)additivity properties. More precisely, those functions $F$ will be defined on a sequence $(G_n)$ of finite graphs with uniform vertex degree bound $d$. Our main goal is to study the normalized limit behaviour, i.e.\@ to determine a class of sequences for which the limit
\[
\overline{F} := \lim_{n\rightarrow \infty} \frac{F(G_n)}{|V(G_n)|}
\]  
exists in the topology of the Banach space. This question is motivated by various mathematical applications: 
\begin{itemize}
\item For subadditive functions (with $\RR$ as underlying Banach space), convergence theorems of the above kind can be used to show the existence of invariants in topological dynamical systems. Those quantities comprise topological entropy and topological mean dimension, see \cite{Gromov-99, LindenstraussW-00} in the situation of amenable groups and \cite{CeccheriniKC-12} for amenable, left cancellative semigroups. 
 \item The notion of almost-additivity for general Banach space-valued functions $F$ has been introduced in \cite{Lenz-02}. Defining $F$ on the set of associated words in a subshift dynamical system, the author characterizes unique ergodicity of the subshift by the normalized Banach space convergence of $F$. A related result for Banach space valued functions on tilings of Euclidean spaces can also be found in \cite{LenzS-05}.
Further, it turns out that this class of almost-additive function is also useful for spectral issues for operators on discrete structures. In detail, it allows for the (uniform) approximation of the spectral distribution function ({\em integrated density of states}, IDS) via finite volume analogues, cf.\@ \cite{LenzS-05, Elek-08}. In the context of abstract amenable groups, this has for instance been done in \cite{KlassertLS-03, Veselic-05b, LenzSV-10, PogorzelskiS-12, Pogorzelski-12}. Related results can also be found in \cite{MuellerR-07, MuellerR-13}, where the authors prove almost-additive convergence theorems that characterize unique ergodicity of Delone dynamical systems. Those assertions lead to finite-dimensional spectral approximations as well.   
  \end{itemize} 

However, all of those convergence theorems deal with F{\o}lner or van Hove sequences in amenable groups or semigroups. 
In the present paper, we turn to a very general class of graph sequences including F{\o}lner sequences in finitely generated groups. Proving convergence theorems along those objects, we significantly extend the geometric framework of the above results. We do so by identifying subadditive and almost-additive functions as mappings defined on graphs satisfying a continuity property with respect to a particular distance function introduced in \cite{Elek-08}. 
Precisely, we determine sufficient and necessary conditions for naturally occurring graph sequences to be Cauchy in this distance function. 
Our proof is based on the groundwork presented in \cite{Elek-12}. Another independent proof 
via algorithmic techniques has recently given by {\sc Newman} and {\sc Sohler}, cf.\@
\cite{NewmanS-13}, Theorem 3.1.  

 This leads to our first key result given in Theorem~\ref{thm:main}. As an immediate consequence, we are able to derive strong Banach space-valued convergence results given in Corollary~\ref{cor:weak} (almost-additive functions) and Theorem~\ref{thm:subadditive} (subadditive functions) for so-called {\em hyperfinite} (and weakly convergent, see below) graph sequences as introduced in \cite{Elek-07}. Roughly speaking, hyperfiniteness of a family of graphs means that for all elements in the family, the proportion of edges that need to be deleted in order to obtain (edge-)disjoint components with a uniform (vertex) size bound, is small uniformly over the family. As a main example, every F{\o}lner sequence in a finitely generated, amenable group is hyperfinite. A more general conjecture of {\sc Elek} (Conjecture~1 in \cite{Elek-08}) refers to bounded vertex degree graphs with edges and vertices labeled by finitely many colours. So far, we only have been able to cover the geometric situation ignoring the colours. Thus, within the uncoloured situation, Corollary~\ref{cor:weak} brings out a unification of all results that concern amenable Cayley graphs, see e.g.\@ \cite{LenzMV-08, LenzSV-10, PogorzelskiS-12}. In addition to that, the graphs modeling quasicrystals in \cite{LenzS-05, LenzV-09} belong just as well to the class to which our convergence theorem applies. Theorem~\ref{thm:subadditive} can be interpreted as a version
of the 'Ornstein-Weiss' lemma which is a crucial tool in the world of group dynamics,
see e.g.\@ \cite{Gromov-99, LindenstraussW-00, Krieger-10}. Our result extends an assertion by {\sc Gromov} in \cite{Gromov-99}. Accepting a monotonicity assumption on the function under consideration, it can be also be understood as an upgrading of \cite{Krieger-10, CeccheriniKC-12} in the situation of sequences (in the latter papers, the authors are also able to deal with F{\o}lner {\em nets}).  \\  
Before we explain the content of this paper in more detail, let us briefly overview the classes of graph sequences which are crucial in our elaborations. \\
The original notion of convergent graph sequences $(G_n)$ was introduced in the seminal work of {\sc Benjamini} and {\sc Schramm}, cf.\@ \cite{BenjaminiS-01}. It means the existence of the occurrence frequencies for arbitrary geometric patterns as $n\rightarrow\infty$. In the present work, we will refer to those sequences as {\em weakly convergent}. In general, the limit is given by a so-called {\em graphing} which is a Borel probability space consisting of countable, uniformly vertex-degree bounded graphs, see e.g. \cite{Elek-071}. In particular situations, graph sequences might also converge to one deterministic, countable graph (almost-surely). This is for instance the case for sofic approximations of finitely generated groups, cf.\@ \cite{Weiss-00, Elek-06, Elek-08, AbertGV-13}.  

The notion of weakly convergent graph sequences is widely used in the context of various mathematical problems. Let us mention some examples: 
\begin{itemize}
\item In \cite{BenjaminiS-01}, the authors prove that in the situation of planar graph sequences, almost all elements in the graphing are recurrent. 
\item Further, graph sequences turn out to be useful in order to show that certain geometric quantities can be obtained as a limit of their finite-dimensional restrictions. A very prominent question in this field is given by the so-called L\"uck conjecture concerning the approximation of the von Neumann dimensions on finitely generated groups. The L\"uck conjecture was proven to be true for very large classes of geometries, see e.g.\@ \cite{Lueck-94, DodziukLMSY-03, Elek-06, Thom-08, AbertTV-13}. 

\item Moreover, the theory of weakly convergent graph sequences has been used before to draw connections between the geometry and the spectral properties of discrete structures. In this context, in \cite{AbertGV-13}, the authors characterize amenability of subgroups in a finitely generated group via the limit behaviour of the spectral radii evolving from a convergent sequence of Schreier graphs.

\item It is also well-known that the distribution function of the spectrum of selfadjoint operators on graphs can be approximated by finite-volume restrictions. In \cite{Elek-08}, the author defines a class of operator sequences defined over the elements of a weakly convergent graph sequence. It turns out that the normalized spectral distributions converge in the sense of weak convergence of measures. For the approximations of sofic, finitely generated groups, a related result can be found in \cite{SchwarzenbergerS-13}. In fact, a group is sofic if it contains a weakly convergent graph sequence approximating the local structures of the corresponding Cayley graph.  
\end{itemize}
Considering further kinds of convergence for graph sequences, there is no unified mathematical theory so far. However, for certain applications, it is essential to find refinements of the Benjamini-Schramm convergence, see e.g.\@ \cite{HatamiLS-13}. As it is our goal to prove uniform spectral approximation results for operators on discrete structures, we will have to deal with a notion which strongly takes into account the overall geometry of the graphs under consideration. It turns out that the 'right' criterion relies on a metric introduced by {\sc Elek} in \cite{Elek-08}.
We show that convergence in this metric (so-called {\em strong convergence}) is equivalent to the sequence being weakly convergent {\em and} hyperfinite. As a consequence, we are able to derive the uniform convergence of the spectral distributions of finite range operators on weakly convergent, hyperfinite graph sequences. As mentioned before, this widely opens the range of geometries that can be considered. However, in the geometrically more restricted situation of amenable groups (F{\o}lner sequences), one may deal with colourings allowing for more general dynamical systems, see e.g.\@ \cite{LenzSV-10, PogorzelskiS-12}. \\ 
Let us briefly summarize the organization of this paper. In Section~\ref{sec:prelim}, we give the necessary preliminaries and we introduce the notions of weak and strong convergence for graph sequences. Further, we define hyperfiniteness for graph sequences and we provide some examples. Next, we prove in Section~\ref{sec:hyper}
that for a weakly convergent sequence, hyperfiniteness is a necessary {\em and} sufficient condition for its strong convergence, cf.\@ Theorem~\ref{thm:main}. The proof relies on the so-called Equipartition Theorem of {\sc Elek}, cf.\@ \cite{Elek-12}. The following two Sections~\ref{sec:BAP} and~\ref{sec:subadditive} are devoted to the proofs of the announced Banach space valued convergence theorems. In Corollary~\ref{cor:weak}, we show the validity of the above limit relation for almost-additive functions. The corresponding result for subadditive functions can be found in Theorem~\ref{thm:subadditive}. Moreover, we provide an operator theoretic framework for an application of the almost-additive convergence thorem, in Section~\ref{sec:IDS}. In detail, just as well as in \cite{Elek-08}, we obtain the uniform convergence of the spectral distribution functions for bounded, self-adjoint, finite hopping range operators towards an abstract notion of integrated density of states (IDS) of the limit object (Theorem~\ref{thm:IDS}).

\section{Graph convergence} \label{sec:prelim}

For the main purposes of this work, we introduce in this section the notions of graph isomorphism classes and explain the concept of weakly convergent and strongly convergent graph sequences. Before doing so, we clarify what we mean by graphs and we give the basic definitions related to these objects. \\

In the following, we will always deal with countable (in most cases in fact finite) graphs. These latter elements can be represented as a pair $G:=(V(G), E(G))$, where $V(G)$ is a countable set called {\em vertex set} and $E(G) \subseteq \{\{x,y \}\,|\, x,y \in V(G), x \neq y\}$ is the so-called {\em edge set} describing abstract connections between pairs of vertices. If $\{x,y\} \in E$, we say that $y$ is a neighbour of $x$ and vice versa. The number of neighbours of a vertex $x \in V(G)$ is called the {\em vertex degree} in $x$. Note that each countable graph can be endowed in a natural way with the {\em path distance metric} $d_G$. Thus, if for $x,y \in V(G)$, there is a finite number of distinct vertices $z_j \in V(G)$, $j=1,\dots n$ such that $\{x,z_1\}, \{z_j, z_{j+1}\}, \{z_n,y\} \in E(G)$, $(j=1,\dots, n-1)$, we say that there is a path of length $(n+1)$ in $G$ joining $x$ and $y$ (and vice versa). In this case, the value $d_G(x,y)$ is given by the minimum of the lengths of paths joining $x$ with $y$. As a convention, we set $d_G(x,y)= \infty$ if there is no such path. Next, we define first the notion of subgraphs which will become essential in the considerations below. 

\begin{Definition}[subgraphs] \label{defi:subgraphs}
Let $G$ be a countable graph with vertex set $V(G)$ and edge set $E(G)$. We call $H=(V(H), E(H))$ a subgraph of $G$ if
\begin{itemize}
\item $V(H) \subseteq V(G)$ and
\item $E(H) \subseteq \{ \{x,y\} \in E(G)\,|\,x,y \in V(H) \}$.
\end{itemize}
\end{Definition}

A special class of subgraphs of a graph $G$ is given by graphs induced by a subset of the vertex set $V(G)$.

\begin{Definition}[induced graph]
Let $G$ be a countable graph with vertex set $V(G)$ and edge set $E(G)$. Then, for each set $T \subseteq V(G)$, we shall define {\em the graph on $T$ induced by $G$} as $G(T):= (V_T, E_T)$, where 
\begin{eqnarray*}
V_T &:=& T, \\
E_T &:=& \{\{x,y\} \in E(G)\,|\, x,y \in T\}.
\end{eqnarray*}
\end{Definition}

Note that induced subgraphs can be interpreted as subgraphs in the sense of Definition \ref{defi:subgraphs} with maximal possible number of edges (for a given vertex set). \\
 
Following a concept introduced by {\sc Elek}, see \cite{Elek-12}, we now turn to convergent graph sequences. To do so, let $d \in \NN$ and denote by $\cG_d$ a set of finite graphs of vertex degree bounded by $d$. We assume further that the graphs in $\cA_d$ do not have loops and that there are no multiple edges. We say that $G \in \cG_d$ is a rooted graph of radius bounded by $r \in \NN$ if there is one distinguished vertex $x := \operatorname{rt}(G) \in V(G)$ (the {\em root} of $G$) such that $d_G(x,y) \leq r$ for every $y \in V(G)$, where $d_G$ is the canonical shortest path metric in $G$, see above. In this sense we say that two graphs $G_1$ and $G_2$ are rooted isomorphic if there is a bijection $\varphi:V(G_1) \rightarrow V(G_2)$ mapping root to root and where $\{\varphi(x), \varphi(y)\} \in E(G_2)$ if and only if $\{x,y\} \in E(G_1)$. We will write $G_1 \sim G_2$ if $G_1$ and $G_2$ are rooted isomorphic. \\
We denote by $\cA_d$ the set of all rooted isomorphism classes of graphs of vertex degree bounded by $d$ and we define the radius of the classes $\alpha \in \cA_d$ 
\[
\operatorname{rad}(\alpha) := \min\{ t \in \NN\,|\, d_G(\operatorname{rt}(\alpha), y) \leq t, \, y \in V(\alpha)\}.
\] 
For $r \in \NN$, we write $\cA_d^r$ for the elements in $\cA_d$ with radius exactly $r$. Note that it follows from the fact that $\cA_d = \bigcup_{r=1}^{\infty} \cA_d^r$ that $\cA_d$ is a countable set. 

Now for $G \in \cG_d$, $r \in \NN$ and $\alpha \in \cA_d^r$, we set
\[
T_r(G,\alpha) := \{x \in V(G)\,|\, B_r^G(x) \sim \alpha \},
\]    
where $B_r^G(x):=G(B_r(x))$ is the induced subgraph of $G$ with vertex set $B_r(x):= \{y \in V(G)\,|\, d_G(x,y) \leq r \}$, i.e.\@ the $r$-ball centered $x \in V(G)$. Further, for each $\alpha \in \cA_d$, the number
\[
p(G, \alpha) := \frac{|T_{\operatorname{rad}(\alpha)}(G,\alpha)|}{|V(G)|}
\]
can be interpreted as the probability that the $\operatorname{rad}(\alpha)$-ball around a random vertex of $G$ is rooted-isomorphic to the element $\alpha$. If we enumerate the elements of the (countable) set $\cA_d$, we get a map
\[
\cL: \cG_d \rightarrow [0,1]^{\NN}: \cL(G) := (p(G, \alpha_i))_{i \in \NN},
\]  
cf. \cite{Elek-12}. Further, we equip $[0,1]^{\NN}$ with a metric $d_{\pi}$ generating the usual product topology. With no loss of generality, we represent $d_{\pi}$ as
\[
d_{\pi}(\cL(G),\cL(H)):= \sum_{k=1}^{\infty} 2^{-k}\, \frac{|p(G, \alpha_k) - p(H, \alpha_k)|}{1 + |p(G, \alpha_k) - p(H, \alpha_k)|}
\] 
for $G, H \in \cG_d$. Note that the map $\cL$ is 'almost' injective in the sense that $\cL(G) = \cL(H)$ implies that there must be a graph $K \in \cG_d$ such that both $G$ and $H$ are disjoint unions of $K$-copies,  cf.\@ \cite{Elek-12}. 

We are now in position to introduce notions of convergence for sequences $(G_n)_{n=1}^{\infty} \subseteq \cG_d$. 

\begin{Definition}[Weak convergence of graphs]
We say that a sequence $(G_n) \subseteq \cG_d$ converges {\em weakly} if for all $r \in \NN$ and every $\alpha \in \cA_d^r$, the limit
\[
p(\alpha) := \lim_{n\rightarrow\infty} p(G_n, \alpha)
\]
exists. So $(G_n)_{n=1}^{\infty}$ is weakly convergent if and only if $(\cL(G_n))_{n=1}^{\infty}$ is convergent pointwise, i.e.\@ it is a Cauchy sequence with respect to the the distance $d_{\pi}$. 
\end{Definition}
In fact, the sequence $(G_n)$ consisting of {\em finite} graphs converges to a {\em probability distribution} of possibly infinite graphs. 
This notion of convergence for graphs has been introduced in the influential work of {\sc Benjamini} and {\sc Schramm}, cf.\@ \cite{BenjaminiS-01}. So we will equivalently refer to weak convergence as {\em Benjamini-Schramm convergence}.\\


By defining a particular metric $\delta_{\rho}$ on $\cG_d$, {\sc Elek} introduced the notion of {\em strongly convergent graph sequences} in \cite{Elek-08}. We briefly explain the construction of $\delta_{\rho}$.
Assume first that the graphs $G,H \in \cG_d$ can be represented on the same vertex set $V$. Let
\[
\delta(G,H) := \frac{|\{v \in V\,|\, S^G(v) \neq S^H(v) \}|}{|V|},
\]
where $S^G(v)$ stands for the $1$-ball (the {\em star}) around $v \in V$ in the graph $G$. In this context, the $\neq$-sign must be understood not only in the sense of rooted isomorphism classes, but it also takes into account the differences in the vertex numberings (one can assume that all vertices in $G$ and $H$ are marked by numbers $\{1,\dots,|V|\}$) of both stars. 
Then, $\delta$ defines a metric on the graphs $G \in \cG_d$ which are defined over the vertex set $V$ (cf.\@ \cite{Elek-08}, Lemma 2.1). As a convention, we set $\delta(G,\epsilon) = 1$ for all $G \in \cG_d$, where $\epsilon$ shall be the empty graph. In a next step, one can come up with a metric which is invariant under permutations of the vertex numbering. In Lemma 2.2 of \cite{Elek-08}, it is shown that for the graphs defined on the same vertex set $V$, the expression
\[
\delta_S(G,H) := \min_{\sigma \in S(V)} \delta(G, H^{\sigma})
\] 
is a well-defined metric, where $S(V)$ is the set of all possible permutations of the vertex numbering in $V$ and $H^{\sigma}$ is the graph obtained by renumbering the vertices, also making sure that the edge relations are preserved. Finally, for graphs $G,H$ which are not necessarily defined on the same vertex set, the so-called {\em geometric distance} 
\[
\delta_{\rho}(G,H) := \inf_{\{q,p \in \NN \,|\, q|V(G)| = p|V(H)|\}} \delta_S(qG, pH)  
\]
has been introduced in \cite{Elek-08}, where the graph $qG$ is represented by $q$ disjoint copies of the graph $G$. By Proposition 2.1 of the same work, $\delta_{\rho}$ defines a metric on the set of isomorphism classes in $\cA_d$. For some finite graph $G \in \cG_d$, we let $\alpha(G)$ be the corresponding isomorphism class in the set $\cA_d$.

\begin{Definition}[Strong convergence of graphs]
We say that a sequence $(G_n) \subseteq \cG_d$ converges {\em strongly}, if $(\alpha(G_n))_{n=1}^{\infty}$ is a Cauchy sequence in the $\delta_{\rho}$-metric.
\end{Definition}  

It is a well-known fact that strong convergence of graphs implies weak convergence, see e.g.\@ \cite{Elek-08}, Proposition 2.2. We will show in the next section that under the additional assumption of {\em hyperfiniteness}, the converse is also true.

\begin{Definition}
A family $\cP \subseteq \cG_d$ is called {\em hyperfinite} or {\em amenable} if for every $\varepsilon > 0$, one can find a number $K_{\varepsilon} \in \NN$ such that we can remove from every $G \in \cP$ less than $\varepsilon\,|E(G)|$ edges such that the resulting graph consists of disjoint components containing at most $K_{\varepsilon}$ vertices. 
\end{Definition}

\begin{Remark}
In the literature, hyperfinite graphs are also referred to as so-called {\em anti-expanders}. 
\end{Remark}

\begin{Examples} \label{exa:hyperfinite}
\begin{itemize}
\item Let $\cP = \{P_n\}_{n=1}^{\infty}$, where $P_n$ is a path of length $n$. Then $\cP$ is hyperfinite. One may e.g.\@ choose $K_{\varepsilon} = 2/\varepsilon$ for $\varepsilon > 0$. 
\item Let $\cP = \{F_n\}_{n=1}^{\infty}$, where $(F_n)_{n=1}^{\infty}$ are the Cayley graphs of the sets in a  F{\o}lner sequence in a finitely generated, amenable group. Using the {\sc Ornstein/Weiss} $\varepsilon$-quasi tiling theory (cf.\@ \cite{OrnsteinW-87, PogorzelskiS-12}), it can be seen that $\cP$ is hyperfinite. 
\item Let $\cP = \{G_n\}_{n=1}^{\infty}$ be a graph sequence of {\em subexponential growth}, i.e.\@ for all $n \in \NN$, for each $x \in V(G_n)$ and every $r \in \NN$, we have $|B^{G_n}_r(x)| \leq f(r)$, where $f:\NN \rightarrow \NN$ is a function which has {\em subexponential growth}. Precisely, this latter property means that for each $\beta > 0$, there exists a number $r_{\beta} > 0$ such that $f(r) \leq \exp(\beta\,r)$ whenever $r \geq r_{\beta}$, cf.\@ \cite{Elek-08}. Then, $\cP$ is a hyperfinite family.  
\end{itemize}
\end{Examples}

In fact, one can verify that strongly convergent graph sequences are hyperfinite. 

\begin{Proposition}[cf.\@ \cite{Elek-08}, Proposition 2.3] \label{prop:stronghyper}
If $(G_n) \subseteq \cG_d$ is strongly convergent, then the family $\cP = \{G_n\}_{n=1}^{\infty}$ is hyperfinite. 
\end{Proposition}


\section{Hyperfiniteness and weak convergence} \label{sec:hyper}

The fundamental insight that two statistically alike graphs with the same number of vertices
are also close in their geometry (i.e.\@ w.r.t.\@ $\delta_{\rho}$ as introduced before) has 
been obtained in \cite{Elek-12} (Theorem~5) and independently in \cite{NewmanS-13} (Theorem~3.1).
We deduce from this that hyperfinite, weakly convergent graph sequences are in fact strongly convergent. This confirms a conjecture of {\sc Elek} in \cite{Elek-08} for the case of graphs with uncoloured vertices and colourless edges. 
This section is devoted to a proof which relies on techniques of 
{\sc Elek's} work, cf.\@ \cite{Elek-12}.  

\begin{Theorem} \label{thm:main}
Let $\cP:= (G_n) \subseteq \cG_d$ be a weakly convergent sequence of graphs which is also hyperfinite. Then, $(G_n)$ is in fact strongly convergent. 
\end{Theorem}

For the proof of Theorem \ref{thm:main}, we need the so-called {\em Equipartition Theorem} stating that statistically similar graphs in a hyperfinite family can be partitioned similarly, cf.\@ \cite{Elek-12}, Theorem~4. 

\begin{Theorem}[{\sc Elek}, Equipartition Theorem, cf.\@ \cite{Elek-12}, Theorem 4] \label{thm:equipart}
Let $\cP \subseteq \cG_d$ be a hyperfinite family. Then, for any $\varepsilon > 0$, there exists $K_{\varepsilon} \in \NN$ with the following property: for any $\beta > 0$, there exists $\delta > 0$ such that if $G \in \cP$ and $H \in \cG_d$ with $d_{\pi}(\cL(G),\cL(H)) \leq \delta$, then there is a way to remove less than $2\varepsilon|E(G)|$ egdes of $G$ and less than $2\varepsilon|E(H)|$ edges of $H$ such that
\begin{itemize}
\item in the remaining graphs $G^{'}$ and $H^{'}$, all connected components have size at most $K_{\varepsilon}$,
\item $\sum_{\alpha \in \cA_d:\, |V(\alpha)| \leq K_{\varepsilon}} \, |c_{\alpha}^{G^{'}}- c_{\alpha}^{H^{'}}| < \beta$,
\end{itemize}
where $C_\alpha^{G^{'}}$ is the set of points that are in a component of $G^{'}$ isomorphic to $\alpha$, and $c_{\alpha}^{G^{'}}:= |C_{\alpha}^{G^{'}}|/|V(G^{'})|$.
\end{Theorem}  


We now prove Theorem \ref{thm:main}. \\

{\bf Proof of Theorem~\ref{thm:main}.}
Let $\varepsilon > 0$. Further, choose $\varepsilon_1:= \varepsilon /(6d)$. For any $n \in \NN$, we remove $\varepsilon_1|E(G_n)|$ edges from $G_n$ such that in the remaining graphs $\{G_n^{'}\}_{n=1}^{\infty}$, all connected components have at most $K_{\varepsilon_1}$ vertices, where the number $K_{\varepsilon_1}$ is chosen according to the Equipartition Theorem~\ref{thm:equipart}. For $n \in \NN$, we call a vertex $x \in G_{n}^{'}$ {\em exceptional} if we removed at least one of the edges incident to $x$. Hence, there can be at most $2d\varepsilon_1|V(G_n)|$ exceptional vertices in $G_n$. We denote by $\alpha_1, \alpha_2, \dots, \alpha_{M_{\varepsilon_1}}$ $(M_{\varepsilon_1} \in \NN)$ the isomorphism classes in $\cA_d$ consisting of at most $K_{\varepsilon_1}$ vertices. So for $n \in \NN$ and $1 \leq i \leq M_{\varepsilon_1}$, we write $\kappa_i^{n} \in \NN$ for the number of connected components in $G_n^{'}$ which are rooted isomorphic to $\alpha_i$. 
We further set $\gamma_i^{n}:= \kappa_i^{n}/|V(G_n)|$ for $n \in \NN$ and $1 \leq i \leq M_{\varepsilon_1}$.
Note that due to the weak convergence, the sequence $(G_n)$ is Cauchy in the $d_{\pi}$-distance and the Equipartition Theorem~\ref{thm:equipart} is applicable. 
Therefore, using that $\kappa_i^{n}:= |C_{\alpha_i}^{G_n^{'}}|/|V(\alpha_i)| \leq |C_{\alpha_i}^{G_n^{'}}|$, where $C_{\alpha_i}^{G_n^{'}}$ is defined as in the Equipartition Theorem,  we find for $\beta:= \varepsilon_1 / ( M_{\varepsilon_1}\,K_{\varepsilon_1})$ some number $L \in \NN$ such that
\[
\left| \gamma_i^{n} - \gamma_i^{m} \right| < \beta = \frac{\varepsilon}{6d\,M_{\varepsilon_1}\,K_{\varepsilon_1}}
\]
for all $n,m \geq L$ and all $1 \leq i \leq M_{\varepsilon_1}$.
Next, using Lemma 2.3 of \cite{Elek-08}, we conclude that indeed, 
\[
d_{\rho}(G_n,G_m) \leq 4d\,\varepsilon_1 + 2\,M_{\varepsilon_1}K_{\varepsilon_1}\beta \leq \varepsilon 
\]
for all $n,m \geq L$. \hfill $\Box$         

\section{A Banach-space valued convergence theorem} \label{sec:BAP}

This section is devoted to the proof of a general convergence result for so-called {\em almost additive functions} on $\cG_d$ which take their values in an arbitrary Banach space. Precisely, we show that for strongly convergent graph sequences $(G_n) \subseteq \cG_d$, the expression $F(G_n)/|V(G_n)|$ converges in the Banach space topology as $n$ tends to infinity, where $F$ belongs to the class of functions under consideration. 
Similar results in the situation of general (not necessarily finitely generated) amenable groups can already be found in \cite{LenzSV-10, PogorzelskiS-12} and \cite{Pogorzelski-12}.\\  

We start with the definition of almost-additive functions on $\cG_d$. Let $X$ be a Banach space equipped with a norm $\|\cdot\|$. 

\begin{Definition} \label{defi:AP}
A mapping
\[
F:\cG_d \rightarrow (X,\|\cdot\|)
\]
is called {\em almost-additive} (on $\cG_d$) if $F(\epsilon) = 0$ (in case the empty graph $\epsilon \in \cG_d$) and if 
there is a constant $D > 0$ depending on $d$ such that if $G,H \in \cG_d$ and $p,q \in \NN$ are such that $p\,|V(G)| = q\,|V(H)|$, then
\[
\left\| p\,F(G) - q\,F(H) \right\| \leq D\,\delta_S(pG, qH)\,p\,|V(G)|.
\]
We call this the {\em almost-additivity} property of $F$.
\end{Definition} 

The following statements are straight forward consequences from Definition \ref{defi:AP} and the considerations of the previous Section \ref{sec:hyper}.

\begin{Theorem} \label{thm:strong}
Assume that $F$ is an almost-additive mapping on $\cG_d$ with values in a Banach space $(X,\|\cdot\|)$. Then, if $(G_n) \subseteq \cG_d$ is a strongly convergent graph sequence, then there is an element $\overline{F} \in X$ such that
\[
\lim_{n\rightarrow \infty} \left\| \frac{F(G_n)}{|V(G_n)|} - \overline{F} \right\| =0. 
\]
\end{Theorem}

{\bf Proof.}
Choose $\varepsilon > 0$ and find $L \in \NN$ such that $d_{\rho}(G_n, G_m) \leq \varepsilon$ whenever $n,m \geq L$. By the definition of $d_{\rho}$, for each choice of integers $n$ and $m$, we find $p, q \in \NN$ with $p |V(G_n)| = q|V(G_m)|$ such that $d_S(pG_n, qG_m) \leq 2\,\varepsilon$. Consequently, we can use the almost-additivity property of $F$ to compute
\begin{eqnarray*}
\left\| \frac{F(G_n)}{|V(G_n)|} - \frac{F(G_m)}{|V(G_m)|} \right\| &=& \left\| \frac{p\,F(G_n)}{p\,|V(G_n)|} - \frac{q\,F(G_m)}{q\,|V(G_m)|} \right\| \\
&\leq& 2D\,\varepsilon
\end{eqnarray*}
for all $n, m \geq L$. It follows that $(F(G_n)/|V(G_n)|)_{n=1}^{\infty}$ is a Cauchy sequence in the Banach space $X$ and hence convergent to some element $\overline{F} \in X$. \hfill $\Box$\\

The following corollary demonstrates that in the situation of translation invariance of the almost-additive function under consideration, the convergence result of \cite{PogorzelskiS-12} along F{\o}lner sequences in a finitely generated, amenable group can be extended to graph sequences. In fact, Definition~\ref{defi:AP} includes periodicity of the function $F$ in the sense that $F(G)=F(H)$ if $G$ is rooted isomorphic to $H$. Introducing graph colourings in \cite{PogorzelskiS-12},   
the authors are able to work with a slightly weaker notion of invariance. However, the corresponding geometric situation is far more restricted as the graph sequences are given as subgraphs of the Cayley graph of the group induced by the elements of the F{\o}lner sequence under consideration. This latter graph sequence is in fact weakly convergent with limit probabilities 0 or 1. Moreover, as we have mentioned in the Examples \ref{exa:hyperfinite} above, the sequence is also hyperfinite. It follows from this that in the translation invariant situation, the setting in \cite{PogorzelskiS-12} is a special case of the next corollary.        

\begin{Corollary} \label{cor:weak}
Assume that $F$ is an almost-additive mapping on $\cG_d$ with values in a Banach space $(X,\|\cdot\|)$. Then, if $(G_n) \subseteq \cG_d$ is a weakly convergent, hyperfinite graph sequence, then there is an element $\overline{F} \in X$ such that
\[
\lim_{n\rightarrow \infty} \left\| \frac{F(G_n)}{|V(G_n)|} - \overline{F} \right\| =0. 
\]
\end{Corollary}

{\bf Proof.}
This follows directly from the previous Theorem \ref{thm:strong} and Theorem \ref{thm:main}.  \hfill $\Box$ \\

\section{Subadditive convergence} \label{sec:subadditive}

Note that Fekete's Lemma is a well known, elementary statement providing essential applications in various mathematical areas. It reads as follows. 

\begin{Lemma} \nonumber
Let $(a_n)_{n \in \NN} \in \mathbb{R}^{\mathbb{N}}$ be a subadditive sequence, i.e.
\[
a_{m+n} \leq a_m + a_n 
\]  
for all $n,m \in \NN$. Then the sequence $(\frac{a_n}{n})$ converges to its infimum
(which might be $-\infty$).
\end{Lemma}

For the sake of applications, one may raise the question whether the integer index set can be replaced by more complicated structures, such as sets or graphs. Quite recently, a corresponding result was proven in \cite{CeccheriniKC-12}, where the authors prove subadditive convergence for amenable, cancellative semigroups in order to apply it to problems concerning the entropy of measure preserving dynamical systems. Including a natural additional condition on the function $h$ under consideration, we prove convergence along hyperfinite graph sequences. Precisely, we need the subadditive function to be non-decreasing with respect to the subgraph ordering. This extends the geometric situation in results obtained by {\sc Lindenstrauss} and {\sc Weiss}, cf.\@ \cite{LindenstraussW-00} as well as by {\sc Gromov}, cf.\@ \cite{Gromov-99} in the context of amenable groups.      
Let us clarify first what we mean by subadditive functions. \\
Assume that $h:\cG_d \rightarrow \RR$ is a {\em subadditive function}, i.e.\@ it satisfies the following properties:
\begin{itemize}
\item there exists a constant $C > 0$ such that $h(G) \leq C\,|V(G)|$ for all $G \in \cG_d$ (boundedness); 
\item if $G \in \cG_d$ is a subgraph of $\tilde{G} \in \cG_d$, then
\[
h(G) \leq h(\tilde{G})  \quad \quad \mbox{(monotonicity);}
\]
\item if $G$ and $G'$ are $\cG_d$-subgraphs of some $\tilde{G} \in \cG_d$ with $V(G)\sqcup V(G') = V(\tilde{G})$, 
\[
h(\tilde{G}) \leq h(G) + h(G') \quad\quad \mbox{(subadditivity);}
\] 
in the {\em edge-disjoint} situation, i.e.\@ if we also have
and $e \in E(\tilde{G})$ if and only if {\em either} $e \in E(G)$ {\em or} \ $e \in E(G')$, then we have in fact
\[
h(\tilde{G}) = h(G) + h(G') \quad\quad \mbox{(special additivity);}
\]
\item it is true that $h(G) = h(G')$ if $G$ is isomorphic to $G'$ in $\cG_d$ (pattern-invariance). 
\end{itemize}

\begin{Remark}
We need the special additivity assumption for the edge-disjoint case as a technical feature of the more general subadditivity poperty. Note that in many settings, this is a pathological situation, e.g.\@ if one considers sequences consisting of connected graphs. However, this property will become relevant in the proofs below, where the strong graph metric $\delta_{\rho}$ forces us to work with $q$-fold ($q \in \NN$) vertex- and edge-disjoint copies of finite graphs $G \in \cG_d$. For those objects, the special additivity assumption and the pattern-invariance condition guarantee that $h(qG) = q\,h(G)$.     
\end{Remark}


\begin{Theorem}\label{thm:subadditive}
Assume that $h:\cG \rightarrow \mathbb{R}$ is a subadditive function. Then, for every weakly convergent, hyperfinite graph sequence $(G_n) \subseteq \cG_d$, there is a number $\lambda \in \RR \cup \{-\infty\}$ such that
\[
\lim_{n \rightarrow \infty} \frac{h(G_n)}{|V(G_n)|} = \lambda. 
\]
\end{Theorem}

{\bf Proof.}
By Theorem~\ref{thm:main}, the sequence $(G_n)$ is in fact strongly convergent. Further, we set 
\[
\lambda := \liminf_{n \rightarrow \infty} \frac{h(G_n)}{|V(G_n)|},
\]
which is an element in the infinite interval $[-\infty, C]$ due to the boundedness of $h$ from above. We denote by $K \subseteq \NN$ an infinite set containing the indices of a subsequence of $(h(G_n)/|V(G_n)|)_n$ that converges to $\lambda$. We fix an arbitrary positive number $\varepsilon > 0$ and we choose $k_0 \in K$ such that for all $n,k \geq k_0$, we have $\delta_{\rho}(G_n, G_k) < \varepsilon$. We fix such $n$ and $k$, where we also make sure that $k \in K$. By definition of $\delta_{\rho}$, there exist integer numbers $q_n, q_k \in \NN$ such that $q_n|V(G_n)| = q_k|V(G_k)|$ and 
\begin{eqnarray} \label{eqn:AA}
\delta_S(q_nG_n, q_kG_k) \leq 2\,\varepsilon.
\end{eqnarray}
Assuming that both $q_nG_n$ and $q_kG_k$ are defined on a common vertex set $V_{n,k}$ we can find a subset $\tilde{V}_{n,k} \subseteq V_{n,k}$ containing those vertices such that the $1$-balls (including vertex numbering) coincide in both graphs. By Inequality (\ref{eqn:AA}), we have $|\tilde{V}_{n,k}| \geq (1-2\varepsilon)|V_{n,k}|$. For this set $\tilde{V}_{n,k}$, we denote the associated induced subgraph in $q_n G_n$ (in $q_k G_k$ respecively) by $\tilde{G}_{n,k}$. Thus, using the special additivity property, as well as subadditivity, boundedness and pattern-invariance of $h$, we obtain
\begin{eqnarray} \label{eqn:compute1}
\frac{h(G_n)}{|V(G_n)|} = \frac{h(q_nG_n)}{q_n|V(G_n)|} \leq \frac{h(\tilde{G}_{n,k})}{|V_{n,k}|} + 2C\,\varepsilon. 
\end{eqnarray}     
Since $h$ is non-decreasing, it follows that $h(\tilde{G}_{n,k}) \leq h(q_kG_k)$. Further, the special additivity property and the pattern-invariance condition of $h$ yield $h(\tilde{G}_{n,k}) \leq q_k h(G_k)$. Continueing Inequality (\ref{eqn:compute1}), we arrive at
\begin{eqnarray*}
\frac{h(G_n)}{|V(G_n)|} \leq \frac{q_k h(G_k)}{|V_{n,k}|} + 2C\,\varepsilon = \frac{h(G_k)}{|V(G_k)|} + 2C\,\varepsilon. 
\end{eqnarray*}   
Note that this latter relation holds true for all large enough $k \in K$ and all large enough $n \in \NN$. Since $\varepsilon > 0$ was arbitrary, we have  
\[
\limsup_{n\rightarrow\infty} \frac{h(G_n)}{|V(G_n)|} \leq \liminf_{n\rightarrow\infty} \frac{h(G_n)}{|V(G_n)|} + 2C\,\varepsilon = \lambda + 2C\,\varepsilon 
\]
for all $\varepsilon > 0$. Hence, sending $\varepsilon\rightarrow 0$ concludes the proof. \hfill $\Box$ \\

\section{Integrated density of states} \label{sec:IDS}

In this section, we demonstrate that the notion of almost-additive, Banach space valued mappings is the right one in order to verify uniform approximation results concerning the integrated density of states for certain operators. More precisely, in the situation of a strongly convergent graph sequence, we show that the spectral distribution functions of the associated operators determined by the local patterns of the graph converge uniformly. As the results of this section have already been proven in \cite{Elek-08}, we do not claim originality, but emphasize the application of the Banach space convergence theorem, Corollary~\ref{cor:weak}. A related IDS result for random Schr\"odinger operators has been proven in \cite{Elek-10}. To complete this section, we add a short discussion concerning an alternative method to derive uniform convergence in general geometric situations. The corresponding result is not new and detailed elaborations can be found in \cite{Thom-08, AbertTV-13}.    

\subsection{The model}

The following setting is deduced from \cite{Elek-08}. 
Suppose that $G$ is a graph with countable vertex set $V(G)$ and with vertex degree bound $d \in \NN$. Assume further that a function
\[
h: V(G) \times V(G) \rightarrow \RR
\]
is given and possesses the following properties. 

\begin{itemize}
\item {\em symmetry}: $h(x,y) = h(y,x)$ for all $x,y \in V(G)$.
\item {\em finite hopping range}: there exists some constant $R> 0$ such that $h(x,y) = 0$ whenever $d_G(x,y) > R$.
\item {\em pattern invariance}: it holds true that for the same constant $R>0$, for each $x \in V(G)$, the function $h_x=h(x,\cdot):V(G) \rightarrow \RR$ is determined by the $R$-neighborhood of $x$. Precisely, this means that if $B^{G}_{R}(x) \sim^{\varphi} \alpha$ for some $\alpha \in \cA_d^{R}$, then there is a map $h_{\alpha}:V(\alpha) \rightarrow \RR$ with
\[
h_x(y) = h(x,y) = h_{\alpha}(\varphi(y))  
\] 
for every $y \in V(B^{G}_{R}(x))$.

\end{itemize}

In this situation, we call the number $R$ the {\em overall range} of the function $h$.
Note that the pattern invariance condition makes sure that $h$ is a translation invariant function. 
Precisely, if there is another $x' \in V(G)$ such that $B_R^{G}(x') \sim B_R^{G}(x)$ with associated isomorphism $\varphi_{x,x'}:V(B_R^{G}(x)) \rightarrow V(B_R^{G}(x'))$, $\varphi_{x,x'}(x) = x'$, then
\[
h(x,y) = h(x', \varphi_{x,x'}(y))
\] 
for all $y \in V(B^{G}_R(x))$. Hence, indeed the function $h$ only depends on the local patterns of the graph $G$.

 In the following, we will refer to functions of the kind as $h$ as {\em admissible kernel functions}. 

\begin{Definition} \label{defi:CO}
Let $(G_n) \subseteq \cG_d$ be a weakly convergent graph sequence. Then, we call $\cH:=(H_n)$ a corresponding {\em weakly convergent operator sequence} if 
\begin{itemize}
\item there is a sequence of admissible kernel functions $(h_n)$, $h_n:V(G_n) \times V(G_n) \rightarrow \RR$, where the parameter $R$, as well as the values of $h_{\alpha}$ $(\alpha \in \cA_d^{R})$ for the $h_n$ {\em do not} change with $n \in \NN$,
\item for all $n \in \NN$, the operator $H_n:\ell^2(V(G_n)) \rightarrow \ell^2(V(G_n))$ is given by
\[
(H_nu)(x):= \sum_{y \in B^{G_n}_R(x)}h_n(x,y)\,u(y). 
\]
\end{itemize} 
\end{Definition} 

\begin{Example}
A simple example for a weakly convergent operator sequence are the discrete Laplacians $\Delta := (\Delta_n)$,  
\[
\Delta_n: \ell^2(F_n) \rightarrow \ell^2(F_n): (\Delta_nu)(x) := \sum_{y\sim x} \left(u(y) - u(x)\right), 
\]
where $(F_n)$ is a sequence of Cayley graphs associated with a F{\o}lner sequence in a finitely generated, amenable group $\Gamma$, cf.\@ \cite{Elek-08} and \cite{PogorzelskiS-12}. 
\end{Example}

\begin{Definition}
We call $(H_n)$ a {\em null operator sequence} if $(h_n)$ is a weakly convergent operator sequence such that for all $\alpha \in \cA_d^{R}$ with $p(\alpha) \neq 0$, we have $h_{\alpha}= 0$.
\end{Definition}

Further, it is well-known that for a fixed, weakly convergent graph sequence $\cG:= (G_n)$, the set $\mathcal{O}_{\cG}$ of corresponding convergent operator sequences possess an algebra structure, i.e.\@ if $H^1$ and $H^2$ are in $\mathcal{O}_{\cG}$, then $H^1 + H^2$ and $H^1 \cdot H^2$ are as well. For a more detailed discussion, the interested reader may e.g.\@ refer to \cite{Elek-08}. In the following, we will write $\cN_{\cG}$ for the set of all nulloperator sequences associated with $\cG$.  
The following proposition is an easy observation which can be proven in the same manner as in \cite{Elek-08}. 

\begin{Proposition}[\cite{Elek-08}, Proposition 3.1]
Let $\cH = (H_n)$ be a convergent operator sequence on a weakly convergent graph sequence $\cG = (G_n)$. Then, the limit
\[
\lim_{n \rightarrow \infty} \frac{1}{|V(G_n)|} \, \sum_{x \in V(G_n)} h_n(x,x)
\] 
exists and is equal to
\[
\operatorname{tr}_{\cG}(\cH):= \sum_{\alpha \in \cA_d^{R}} p(\alpha)\, h_{\alpha}(\operatorname{rt}(\alpha)),
\]
where the parameters $R$, $h_{\alpha}$ are as in Definition \ref{defi:CO}.
\end{Proposition}


We will call $\operatorname{tr}_{\cG}(\cH)$ the {\em trace} of $\cH$ in the algebra $\cO_{\cG}$. The next proposition shows that this latter notion contains valuable information on the structure of the quotient space $\mathcal{Q}_{\cG}:= \cO_{\cG} / \cN_{\cG}$, as well as on its closure with respect to a particular (weak) topology.  

\begin{Proposition}[{\sc Elek}, cf.\@ \cite{Elek-08}, Prop.\@ 3.1, Lemmas 3.2, 3.3, 3.4]
Let $\cG:=(G_n)$ be a weakly convergent graph sequence and assume that $\cH^1:= (H^1_n)$ and $\cH^2:= (H^2_n)$ are arbitrary corresponding convergent operator sequences in $\cO_{\cG}$. Then the following statements hold true. 
\begin{enumerate}[(i)]
\item $\operatorname{tr}_{\cG}(\cH^1\cdot \cH^2) = \operatorname{tr}_{\cG}(\cH^2\cdot \cH^1)$.
\item $\operatorname{tr}_{\cG}$ is faithful on $\cQ_{\cG}$, i.e.\@ $\operatorname{tr}_{\cG}(\cH^1 \cdot \cH^1) > 0$ if $\cH^1 \notin \cN_{\cG}$ and $\operatorname{tr}_{\cG}(\cH^1)= 0$ for $\cH^1 \in \cN_{\cG}$. 
\item the algebra $\cQ_{\cG}$ is a pre Hilbert space with inner product
\[
\langle [\cH^1],[\cH^2] \rangle = \operatorname{tr}_{\cG}(\cH^2\cdot \cH^1), 
\]
where $[\cH^i]$ stands for the class of $\cH^i$ in $\cQ_{\cG}$. 
\item the representation $\mathcal{L}_{[\cH^1]}\,[\cH]:= [\cH^1 \cdot \cH]$ $(\cH \in \cO_{\cG})$ for $[\cH^1] \in \cQ_{\cG}$ is a bounded (multiplication) operator on the quotient space $\cQ_{\cG}$. 
\item The weak closure $\overline{\cQ_{\cG}}$ of $\cQ_{\cG}$ is a von-Neumann algebra. Consequently, the trace
\[
\operatorname{tr}_{\cG}([\cH^1]) = \langle \mathcal{L}_{[\cH^1]}\,\one, \one \rangle
\] 
extends to $\overline{\cQ_{\cG}}$ as a ultraweakly  continuous, faithful trace.   
\item $\mathcal{L}_{[\cH^1]}$ is a bounded, self-adjoint operator on $\overline{\cQ_{\cG}}$ with $\sup_{n \in \NN} \|H^1_n\| \leq K(\cH^1)$ for some constant $K(\cH^1) > 0$.  
\end{enumerate} 
\end{Proposition}


For convenience, we simply write $[\cH]$ instead of $\mathcal{L}_{[\cH]}$ $([\cH] \in \overline{\cQ_{\cG}})$ and we call $[\cH]:=(H_n)$ the limit operator (which actually is defined on $\overline{\cQ_{\cG}}$) of the sequence $\cH \in \cO_{\cG}$. We have seen in the above proposition that $[\cH]$ is a bounded, self-adjoint operator. It follows from the spectral theorem that there must be a decomposition
\[
[\cH] = \int_{\RR} \lambda\,dE^{[\cH]}_{\lambda},
\]
where $E^{[\cH]}_{\lambda} \in \overline{\cQ_{\cG}}$ is the spectral projection $\one_{]-\infty,\lambda]}([\cH])$ for $\lambda \in \RR$. Further, the {\em integrated density of states} of $[\cH]$ can be defined by
\begin{eqnarray*}
N_{[\cH]}(\lambda):= \operatorname{tr}_{\cG}(E_{\lambda}^{[\cH]}). 
\end{eqnarray*}

Given a convergent operator sequence $[\cH]=(H_n)$, then for each $n \in \NN$, we set 
\[
n_{H_n}(\lambda) := \#\{E \leq \lambda \,|\, E \mbox{ is an eigenvalue of } H_n  \} \quad (\lambda \in \RR),
\]
as the cumulative eigenvalue counting function for the (matrix) operator $H_n$. The symbol $\#$ shall express that the eigenvalues are counted {\em with} multiplicities. The corresponding spectral distribution functions are just the cumulative eigenvalue distributions (normalized eigenvalue counting functions) given by
\[
N_{H_n}(\lambda) := \frac{n_{H_n}(\lambda)}{|V(G_n)|}, \quad n \in \NN, \quad \lambda \in \RR.
\]  
With these notions at hand, it is not hard to show the weak convergence of the eigenvalue counting functions towards the integrated density of states.   

\begin{Theorem}[c.f.\@ e.g.\@ \cite{Elek-08}, Theorem 1] \label{thm:IDSweak}
Let $\cH=(H_n)$ be a convergent operator sequence on a weakly convergent graph sequence $\cG:= (G_n)$. Then, for any continuity point $\lambda \in \RR$, we have 
\[
\lim_{n\rightarrow\infty} N_{H_n}(\lambda) = N_{[\cH]}(\lambda)
\]
\end{Theorem}

{\bf Proof.}
See e.g.\@ \cite{Elek-08}, Theorem 1. \hfill $\Box$ \\

\subsection{Uniform approximation of the IDS}

We will now outline how the uniform convergence of the IDS can be obtained by using the Banach space convergence theorem, Corollary~\ref{cor:weak}. 




\begin{Theorem}\label{thm:IDS}
Let $\cH=(H_n)$ be a convergent operator sequence on a weakly convergent, {\em hyperfinite} graph sequence $\cG:= (G_n)$. Then, we have in fact uniform convergence of the spectral distributions, i.e.\@ 
\[
\lim_{n\rightarrow\infty} \|N_{H_n}(\cdot) - N_{[\cH]}(\cdot) \|_{\infty} = 0.
\]
\end{Theorem}

\begin{Remark}
The validity of Theorem \ref{thm:IDS} has already been verified in \cite{Elek-08} (cf.\@ Proposition~3.2). However, we would like to emphasize that it is a different approach is to represent the functions $N_{H_n}(\cdot)$ as abstract, almost-additive mappings on  $\cG_d = \{G_n\}$, where the Banach space under consideration is $(\cC_{rb}(\RR), \|\cdot\|_{\infty})$, i.e.\@ the collection of right-continuous functions endowed with supremum norm. Thus, using Corollary \ref{cor:weak} for the proof of uniform convergence towards the integrated density of states, we also show that the notion of almost-additive mappings is an appropriate abstract description for this spectral issue. 
\end{Remark}

The validity of Theorem \ref{thm:IDS} follows essentially from the following proposition. 

\begin{Proposition} \label{prop:AAdd}
Let $\cH=(H_n)$ be a convergent operator sequence on a weakly convergent graph sequence $\cG:= (G_n)$. Then for $\cG_d:=\{G_n\}$, the mapping 
\[
F:\, \mathcal{G}_d \rightarrow \cC_{rb}(\RR) \quad F(G_n) := n_{H_n}(\cdot)
\]
is almost-additive. 
\end{Proposition}

We will not give a proof of the above proposition as it follows essentially from the definition of the metric $\delta_S$ in combination with a uniform rank estimate for the difference of the eigenvalue counting functions of self-adjoint, finite dimensional operators, see e.g.\@ \cite{LenzSV-10}, Proposition 7.2. Very similar computations can also be found in \cite{LenzV-09}, \cite{PogorzelskiS-12}, \cite{Pogorzelski-12} and \cite{Elek-08}, \cite{Elek-10}. Combining the above proposition with Corollary~\ref{cor:weak}, we have indeed verified the validity of Theorem~\ref{thm:IDS}.\\

It can already be noticed in the group situation that one cannot expect the Banach space convergence of Theorem~\ref{thm:strong} to hold for non-amenable structures. Indeed, it has been shown in \cite{Elek-11} (Proposition~4.1) that sofic approximations of groups are hyperfinite if and only if the group under consideration is amenable. Thus, we infer from Proposition~\ref{prop:stronghyper} that our notion of almost-additive functions will not allow for normalized convergence along approximating sequences in {\em non-amenable} groups. However, as far as the question of pointwise (or uniform) convergence of spectral distribution functions is concerned, one can bypass this problem by extending the algebraic eigenvalue conjecture (see e.g.\@ \cite{DodziukLMSY-03}) to the general geometric case of graph sequences. This issue has been addressed by {\sc Ab{\'e}rt, Thom} and {\sc Vir{\'a}g} in \cite{AbertTV-13}. The method is to use the number theoretic techniques developed in \cite{Thom-08} for the proof of the conjecture in the case of sofic groups. Precisely, if $\nu_n$ and $\nu$ are the positive measures associated with $N_{H_n}$ and $N_{[\cH]}$ respectively, then the following assertion holds true. 

\begin{Theorem}[cf.\@ \cite{AbertTV-13}]
Let $(G_n) \subseteq \cG_d$ be a weakly convergent graph sequence. Assume that $[\cH]:=(H_n)$ is an associated weakly convergent operator sequence with integer coefficients (admissible kernel functions), i.e.\@ for all $\alpha \in \cA_d$, we have that $h_{\alpha}$ takes its values in $\ZZ$. Then,
\begin{eqnarray} \label{eqn:points}
\lim_{n\to\infty} \nu_n(\{\lambda\}) = \nu(\{\lambda\})
\end{eqnarray}
for every $\lambda \in \RR$.
\end{Theorem}

Now well-known arguments can be used to derive from the limit relation (\ref{eqn:points}) the uniform convergence of $(N_{H_n})$, see e.g.\@ \cite{LenzV-09}, Lemma 6.3. Thus, one obtains the uniform existence of the integrated density of states for bounded, integer coefficient operators on all approximable, discrete structures with vertex degree bound.  To the knowledge of the author, this is still an open question for arbitrary coefficients, see the remark following Theorem 4.4 in \cite{Thom-08}. \\
These elaborations show that strong spectral approximation results hold true even if one does not have a Banach space valued convergence theorem at hand. However, one may raise the question whether one can find a class of Banch space valued functions which comprises spectral distributions and which are endowed with suitable conditions for a convergence theorem to hold. This will be a subject of future investigations.

\subsection*{Acknowledgements} 

Many thanks go to Daniel Lenz for raising the issue of Banach space convergence theorems for almost-additive functions along graph sequences. Furthermore, the author acknowledges Daniel Lenz' valuable comments concerning the structuring of the present paper. Moreover, the author thanks Siegfried Beckus for a very fruitful discussion about Theorem~\ref{thm:subadditive}. In addition to this, the author expresses his thanks to Andreas Thom for a very inspiring discussion at Leipzig University, as well as for valuable reference information. The author is highly grateful to Miklos Ab{\'e}rt for useful comments and for drawing the author's attention to the Leipzig Spring School lecture notes. Further, the author would like to thank G{\'a}bor Elek
for an enlightening discussion on hyperfinite graph sequences 
during the workshop on graph limits, groups and stochastic processes, Budapest 2014. 
The author gratefully acknowledges the support from the German National Academic Foundation (Studienstiftung des deutschen Volkes). The results of the present paper are part of the doctoral dissertation of the author.    

\bibliographystyle{abbrv}
\bibliography{PS_lit}

\end{document}